\documentclass[12pt]{article}
\usepackage{amsmath,amssymb}
\newcommand{\qed}{{\unskip\nobreak\hfil\penalty50\hskip2em\vadjust{}
    \nobreak\hfil$\square$\parfillskip=0pt\finalhyphendemerits=0\par}}

\newcommand{\x}{\times}

\newcommand{\bs}{\bigskip}

\newcommand{\C}{\mathbb C}

\newcommand{\D}{\mathbb D}
\newcommand{\e}{\varepsilon}

\renewcommand{\H}{\mathbb H}
\renewcommand{\l}{\lambda}

\newcommand{\mb}{\mbox}

\renewcommand{\O}{\Omega}
\renewcommand{\o}{\omega}

\newcommand{\stt}{\subset}
\renewcommand{\t}{\tau}

\newcommand{{\z}}{\Bbb Z}

\newcommand{\su}{subordination }
\newcommand{\an}{analytic }
\newcommand{\VN}{von Neumann }
\renewcommand{\sb}{subalgebra}

\begin{document}
\setlength{\baselineskip}{16pt}

\begin{center}
{\large\bf Analytic Subordination Consequences
of Free Markovianity}

\bigskip\bigskip
{\sc Dan Voiculescu}\footnote
{Supported in part by National Science Foundation
grant NSF--9500308.}

\bs
Department of Mathematics\\ University of California\\
Berkeley, California 94720--3840

\bs\bs\bs{\bf 1. Introduction}
\end{center}

In [4], under some easy to remove genericity assumptions,
we proved the analytic subordination of the Cauchy transforms
of distributions $G_{\mu_{X+Y}}$ and $G_{\mu_X}$ for a pair of
freely independent self-adjoint random variables $X,Y$.
We used this to obtain inequalities among $p$-norms of
densities of distributions, free entropies and Riesz energies.
This was followed by P.Biane's discovery [1] that
\su extends, roughly speaking, to the resolvents,
i.e.~is an operator-valued occurrence. He also showed that this implies a
noncommutative Markov-transitions property for free increment processes
and that there are similar \su results for multiplicative processes.

The \an function approach in [4] and the combinatorial one in
[1] did not shed much light on why \an \su appears
in this context. In [6] we found a simple explanation of this
phenomenon based on the coalgebra structure associated with the
free difference quotient derivation. In the additive case this also
led to a far reaching generalization of the \an \su result to
the $B$-valued case, i.e.~to free independence with amalgamation
over an algebra $B$.

Here, based on simple operator-valued \an function considerations,
we build further on the result in [6]. We derive more general
\an \su results for a freely Markovian triple $A,B,C$.
This also includes the $B$-valued extension of the result for
multiplication of free unitary variables.

\newpage\begin{center}{\bf 2. Notation and Preliminaries}
\end{center}

Throughout, $(M,\tau)$ will denote a tracial
$W^*$-probability space, i.e.~a \VN algebra endowed with a
normal faithful trace state. A triple of \VN \sb s $A,B,C$
contained in $M$ and containing ${\C}1$ is {\it freely
Markovian} (see [5] or [7]) if $A$ and $C$ are $B$-free
in $(M,E_B)$ where $E_B$ is the canonical conditional expectation
of $M$ onto $B$ (see [7] or [8]). If $A,B$ are \VN \sb s of $M$
containing 1, we shall denote by $A\vee B$ or
$W^*(A\cup B)$ the \VN \sb \ generated by $A\cup B$.
If $X=X^*\in M$ and $1\in B\subset M$ is a $*$-\sb \ we shall
denote by $B[X]$ the $*$-\sb \ generated by $B$ and $X$ and we
shall denote by $E_{B[X]}$ the conditional expectation onto the
\VN \sb \ $W^*(B[X])$. For several selfadjoint elements $X,Y,Z,\dots$
(noncommuting) we shall also use notations
$B[X,Y,Z]$, $E_{B[X,Y,Z]}$, etc.

If $A$ is a unital $C^*$-algebra we denote by
${\H}_+(A)=\{T\in A\mid {\mb{Im }} T\geq \e 1$ for some
$\e>0\}$ the {\it upper half-plane of} $A$ and by
${\H}_-(A)=-{\H}_+(A)$ the {\it lower half-plane.}
Note that
$T\in{\H}_+(A)\Leftrightarrow T^{-1}\in {\H}_-(A)$,
in particular elements in ${\H}_{\pm}(A)$ are invertible.
  Also, $A_{sa}$ will denote the selfadjoint elements in $A$.

The open unit ball will be denoted ${\D}(A)=\{T\in A\mid \|T\|< 1\}$. For
balls of radius $R$ we shall also use the notations
${\D}_R(A)=R{\D}(A)$.

\bs\bs\begin{center}{\bf 3. Analytic Subordination}
\end{center}

We begin by recalling the result
([6] Theorem 3.8) which serves as our starting point.

\bs\noindent{\bf 3.1 \ Theorem.} {\it Let $1\in B\stt M$ be a $W^*$-\sb \ and
let
$X\!=\!X^*\!\in\! M$, $Y\!=\!Y^*\!\in\! M$. Assume $X,Y$ are $B$-free
in $(M,E_B)$. Then there is a holomorphic map\\
$F:{\H}_+(B)\to {\H}_+(B)$ such that}
\[
E_{B[X]}((X+Y)-b)^{-1} \ = \ (X-F(b))^{-1}
\]
{\it if} \ $b\in {\H}_+(B)$.

\bs\noindent{\bf 3.2 \ Proposition.} {\it Let $1\in B\stt M$ and
$1\in A\stt M$ be $W^*$-\sb s and let\linebreak
$X=X^*\in M$. Assume $A$ and $X$ are $B$-free
in $(M,E_B)$. Then there is a holomorphic map
$F:{\H}_+(A)\to {\H}_+(B)$ such that}
\[
E_{B[X]}(a-X)^{-1} = (F(a)-X)^{-1}
\]

\bs
{\bf Proof.} It is easily seen that
\[
F_1(a)=(E_{B[X]} (a-X)^{-1})^{-1}+X
\]
is a well-defined holomorphic map
\[
{\H}_+(A)\to {\H}_+(W^*(B[X])) \ .
\]
We must show that $F_1({\H}_+(A))\stt B$. By Theorem 3.1,
\[
F_1(i\e I+A_{sa})\stt {\H}_+(B)
\]
if $\e >0$. Indeed, if $-Y\in A_{sa}$ then Proposition 3.2 gives
\[
E_{B[X]} (i\e I-Y-X)^{-1} = (b-X)^{-1}
\]
for some $b\in{\H}_+(B)$ which means $F_1(i\e I-Y)=b$.

Since $\e I+A_{sa}$ is a uniqueness set for holomorphic functions
on ${\H}_+(A)$, the inclusion $F_1(i\e I+A_{sa})\stt B$
implies $F_1({\H}_+(A))\stt B$ (use functionals $f\in M_*$,
$f|B=0$ to transform the given inclusion into
$(f\circ F_1)\mid i\e I+A_{sa}=0)$. \qed

\bs\noindent{\bf 3.3 \ Proposition.}
{\it Let $A,B,C$ be a freely Markovian triple of
$W^*$-\sb s in $(M, \t)$. Then there is a holomorphic map}
\[
F:{\H}_+(A)\x {\H}_+(C) \to B
\]
{\it such that}
\[
(E_{A\vee B}(a+c)^{-1})^{-1} = a+F(a,c)
\]
{\it if} \ $a\in {\H}_+(A)$, \ $c\in {\H}_+(C)$.

\bs{\bf Proof.} Clearly, we have $a+c\in {\H}_+(M)$, \
$(a+c)^{-1}\in {\H}_-(M)$, \
$E_{A\vee B}(a+c)^{-1}\in {\H}_-(M)$,
\[
(E_{A\vee B}(a+c)^{-1})^{-1}\in {\H}_+(M) \ .
\]
Hence
\[
F_1: {\H}_+(A)\x {\H}_+(C)\to A\vee B\vee C
\]
given by
\[
F_1(a,c)=(E_{A\vee B}(a+c)^{-1})^{-1}-a
\]
is well-defined and holomorphic. Thus, the proof reduces to
showing that\newline $F_1({\H}_+(A)\x {\H}_+(C))\stt B$.
Since $(i\e I+A_{sa})\x (i\e I+C_{sa})$, where $\e>0$ is a
set of uniqueness for holomorphic maps on ${\H}_+(A)\x {\H}_+(C)$,
it suffices to show that
\[
(E_{A\vee B}((X+Y)+2\e i)^{-1})^{-1}-X\in B
\]
when $X\in A_{sa}$, \ $Y\in C_{sa}$. Since
$(X+Y+2\e i)^{-1}\in W^*(B[X])\vee C$ and $A$ and $C$
are $B$-free in $(M,E_B)$ we have
\[
E_{A\vee B}(X+Y+2\e i)^{-1} =
E_{B[X]}(X+Y+2\e i)^{-1} \ .
\]
By Theorem 3.1 we have
\[
E_{B[X]}(X+Y+2\e i)^{-1} = (X+b)^{-1}
\]
for some $b\in {\H}_+(B)$. In particular
\[
(E_{A\vee B}(X+Y+2\e i)^{-1})^{-1} -X\in B \ .
\] \qed

\bs\noindent{\bf 3.4 \ Lemma.} {\it
If $x\in A$, where $A$ is a unital $C^*$-algebra,
the following are equivalent:}

(i) \ $\|x\| <1$

(ii) \ $1-x$ {\it is invertible and $2{\rm{Re}}(1-x)^{-1}\geq (1+\e)$
for some} $\e>0$.

\bs
{\bf Proof.} This fact, not new, is for instance an immediate
consequence of unitary dilation [3] and of the corresponding fact
when $A={\C}$. For the reader's convenience, here is a
direct proof.

Both (i) and (ii) imply $1-x$ is invertible, in which case we have:
\[
(1-x)^{-1} + (1-x^*)^{-1} =
1+(1-x)^{-1}(1-xx^*)^{-1} (1-x^*)^{-1} \ .
\]
To prove (i)\,$\Rightarrow$\,(ii) remark that
$1-xx^*\geq (1-\|x\|^2)1$ and that
\[
\|(1-x^*)^{-1}\xi\|\geq \|1-x^*\|^{-1}
\|\xi\|\geq 2^{-1}\|\xi\|
\]
when $A$ acts on a Hilbert space ${\mathcal H}$ and
$\xi\in{\mathcal H}$. This gives
$(1-x)^{-1}(1-xx^*) (1-x^*)^{-1}\geq 4^{-1}(1-\|x\|^2)1$.

To show (ii)\,$\Rightarrow$\,(i) note that (ii) implies
$(1-x)^{-1}(1-xx^*)(1-x^*)^{-1}\geq \e 1$ so that
$1-xx^*\geq 0$ and $1-xx^*$ is invertible, which gives (i). \qed

\bs\noindent{\bf 3.5 \ Proposition.}
{\it Let $A,B,C$ be a freely Markovian triple
of
$W^*$-algebras in $(M,\t)$ and let
$\O=\{(a,c)\in A\x C\mid a \ {\rm{invertible }}, \
\|a^{-1}c\|<1\}$. Then there is a holomorphic map
$\Phi:\O\to B$ such that if $(a,c)\in\O$, then
$\|a^{-1}\Phi(a,c)\| <1$ and}
$E_{A\vee B}(a-c)^{-1}=(a-\Phi(a,c))^{-1}$.

\bs{\bf Proof.} Let $(a,c)\in\O$. Then $a-c$ is invertible,
the inverse being $(1-a^{-1}c)^{-1}a^{-1}$. Using Lemma 3.41 we have
\[
E_{A\vee B}(1-a^{-1}c)^{-1}\in \textstyle{\frac 12}-i{\H}_+(M)
\]
and hence
\[
E_{A\vee B}(1-a^{-1}c)^{-1} = (1-x)^{-1}
\]
for some $x\in{\D}(M)$. Moreover there is a holomorphic map
$\Psi:\O\to {\D}(M)$ such that
$\Psi(a,c)=1-(E_{A\vee B}(1-a^{-1}c)^{-1})^{-1}$.

Let $\Phi_1:\O\to M$ be defined by
$\Phi_1(a,c)=a\Psi(a,c)$.  Clearly $\Phi_1$ is  holomorphic,
$\|a^{-1}\Phi_1(a,c)\| < 1$ and
\[
E_{A\vee B}(a-c)^{-1} =
(E_{A\vee B}(1-a^{-1}c)^{-1})a^{-1} =
(1-\Psi(a,c))^{-1}a^{-1} = (a-\Phi_1(a,c))^{-1} \ .
\]
Thus the proof reduces to showing that
$\Phi_1(a,c)\in B$.

The open set $\O$ is connected. Indeed, if $(a,c)\in\O$
then $(\l a,c)\in\O$ for all $\l\in{\C}$, $|\l|\geq 1$.
In particular there is a segment in $\O$ connecting $(a,c)$ and
$(Na,c)$ where $N\geq 1$ is such that
$\|(Na)^{-1}\|\leq (1+\|c\|)^{-1}$. Then
$(Na,tc+(1-t)1)\in\O$ for $0\leq t\leq 1$ and we have a path
from $(a,c)$ to $(Na,1)$. If $Na=u(N|a|)$ is the polar
decomposition, then $N|a|\geq (1+\e)1$ for some $\e>0$
and the elements $(u((1-t)N|a|+t(1+\e)1,1)$,
$0\leq t\leq 1$, are in $\O$ connecting $(Na,1)$ and
$((1+\e)u,1)$. Finally $u=\exp (ih)$ for some
$h=h^*\in A$ and $((1+\e)\exp(ith),1)\in\O$ for $0\leq t\leq 1$
continue our path to $((1+\e)1,1)$, etc.

Since $\O$ is connected, the  set
\[
\o=\{(a,c)\in A\x C\mid \|a-3i1\|<1, \
{\mb{Im}}\,c\leq -{\textstyle{\frac 12}}1, \ \|c\|<1\}
\subset {\H}_+(A)\x{\H}_-(C)
\]
and $\o\stt\O$ is a uniqueness set for holomorphic functions in $\O$.
Moreover, using Proposition 3.3 we have
$\Phi(a,c)=F(a,-c)\in B$ where $(a,c)\in\o$ and
$F$ is the function in Proposition 3.3.

\bs\noindent
{\bf 3.6 \ Theorem.}
{\it Let $1\in B$, $1\in C$ be $W^*$-algebras in $(M,\t)$ and
let
$u$ be a unitary element in $M$. Assume $C$ and $\{u, u^*\}$
are $B$-free in $(M,E_B)$. Then there is a holomorphic map}
\[
G:{\D}(C)\to {\D}(B)
\]
{\it such that}
\[
E_{B[u,u^{-1}]}(u-c)^{-1} = (u-G(c))^{-1}
\]
{\it if} \ $c\in{\D}(B)$.

\bs
{\bf Proof.}
This follows immediately from Proposition 3.5 since
$(u,c)\in\O$ if $c\in{\D}(C)$ and
$\|u^{-1}\Phi(u,c)\|<1\Leftrightarrow \|\Phi(u,c)\|<1$
and we take $G(c)=\Phi(u,c)$. \qed

\bs\bs\begin{center}{\bf References}\end{center}
\begin{description}
\item{[1]} P.Biane, Processes with free increments,
{\it Math. Z.} {\bf 227} (1998), 143--174.
\item{[2]} E.Hille and R.S.Phillips, {\it
Functional Analysis and Semigroups}, vol.~31\\
American Math Society Colloqu. Publ., Rhode Island, 1957.
\item{[3]} B.Sz.-Nagy and C.Foias, {\it Analyse Harmonique des
Op\'erateurs de l'espace de Hilbert}, Masson et Cie,
Paris, 1967.
\item{[4]} D.Voiculescu, D. Voiculescu, The analogues of entropy and 
of Fisher's
information\\ measure in free probability theory, I. {\it Commun. Math.
Phys.} {\bf 155} (1993), 71--92.
\item{[5]} D.Voiculescu,  The analogues of entropy and of Fisher's
information measure in free probability theory, VI:   Liberation and
mutual free information. {\it Adv. Math.}  {\bf 146}
(1999), 101--106.
\item{[6]} D.Voiculescu,  The coalgebra of the free difference
quotient and free probability.\\
{\it International Math. Res. Notices} No.~2 (2000), 79--106.
\item{[7]} \ D.Voiculescu, Lectures on Free Probability Theory,
in {\it Lectures on Probability Theory and Statistics},
Saint-Flour XXVIII (1998) Lecture Notes in Math. vol.~1738, 279--349.
\item{[8]} \ D.Voiculescu, K.Dykema and A.Nica,
{\it Free Random Variables},
CRM Monograph Series vol.~I, American Math Society, Rhode Island,
1992.
\end{description}

\end{document}